\documentclass[11pt,twoside, final]{amsart}
\copyrightinfo{0}{Iranian Mathematical Society}
\usepackage{amsmath,amsthm,amscd,amsfonts,amssymb,enumerate}
\usepackage{graphicx}		
\usepackage{color}
\usepackage[colorlinks]{hyperref}

\theoremstyle{definition}

\newtheorem*{NEW} {Proof of   Theorem 1}
\newtheorem*{prooff} {Proof  of  Theorem 2}

\newtheorem*{thm1}{Theorem       1}
\newtheorem*{thm2}{Theorem       2}
\newtheorem*{prop1}{Proposition  1}
\newtheorem*{pro2}{Proposition  2}

\newtheorem*{lem1}{Lemma   1     }
\newtheorem*{lem2}{Lemma   2     }

\newtheorem*{remark1}{Remark  1}
\newtheorem*{remark2}{Remark  2}
\newtheorem*{remark3}{Remark  3}
\newtheorem*{remark4}{Remark  4}

\theoremstyle{definition}

\theoremstyle{remark}

\numberwithin{equation}{section}

 \commby{Hamid Pezeshk}
\begin{document}


\title[Functional equations on $C^{*}$ algebras]{On a functional equation for symmetric linear operators on $C^{*}$ algebras}

\author[Ali Taghavi]{Ali Taghavi}
\address[Ali Taghavi]{Faculty of Mathematics and Computer Science,  Damghan  University,  Damghan,  Iran.}
\email{taghavi@du.ac.ir}



%
\maketitle
%

\begin{abstract}
Let $A$  be a    $C^{*}$ algebra  and $T: A\rightarrow A$  be  a   linear   map which satisfies the functional  equation $\begin{cases}T(x)T(y)=T^{2}(xy)\\T(x^{*})=T(x)^{*} \end{cases}$ We  prove that under each of the following  conditions, $T$ must be the trivial map  $T(x)=\lambda x$ for  some $\lambda \in \mathbb{R}:$\\
\begin{enumerate}
\item $A$  is a simple  $C^{*}$-algebra.\\

 \item $A$ is  unital with trivial center and has  a faithful trace such that  each

  zero-trace element lies in the  closure of the span of  commutator  elements.\\

\item $A=B(H)$  where H is a  separable Hilbert  space.
\end{enumerate}
\

For a given  field  $F$, we consider a similar functional equation  $$\begin{cases}T(x)T(y)=T^{2}(xy)\\T(x^{tr})=T(x)^{tr} \end{cases}$$ where $T$ is  a  linear  map on $M_{n}(F)$ and "tr" is the transpose operator. We prove that this  functional equation has  trivial solution for all $n\in \mathbb{N}$ if  and only if $F$ is  a  formally real field.
\\

\end{abstract}

\section{\bf Introduction}

Motivated by the classical operator  of  differentiation, we  shall consider the  functional equation  \begin{equation}\label{MA}
\begin{cases}T(x)T(y)=T^{2}(xy)\\T(x^{*})=T(x)^{*} \end{cases}
\end{equation}

where $T$ is  a linear  operator  on  a  $C^{*}$  algebra  $A$. We give some sufficient  conditions under which the equation (\ref{MA}) has  only trivial  solution $T(x)=\lambda x$ for  some $\lambda \in \mathbb{R}$. On the  other   hand, in the particular case $A=M_{n}(\mathbb{C})$,  this  equation gives us the following functional equation on $M_{n}(\mathbb{R})$ or  more  generally on $M_{n}(F)$ where $F$ is  an arbitrary field:  \begin{equation}\label{SE}
\begin{cases}T(x)T(y)=T^{2}(xy)\\T(x^{tr})=T(x)^{tr} \end{cases}
\end{equation} where "$tr$" is the transpose operator.
We observe that for  all $n\in \mathbb{N}$, this functional equation has only trivial scalar solution  $$T(x)=\lambda x,\;\;\text{for some}\; \lambda \in F$$ if and only if $F$ is  a  formally real field, that is $\sum_{i=1}^{n} f_{i}^{2}=0,\;f_{i}\in F,\;n\in \mathbb{N}$ implies that $f_{i}=0,\; \forall i=1,2,\ldots,n$.\\
First we explain that how  we construct our main functional equation (\ref{MA}) from the differentiation  operator:\\

Recall that  a  complex  coalgebra is a complex  vector  space  $C$ with linear  maps  $\Delta :C \rightarrow C\otimes C$ and
 $\varepsilon :C \rightarrow \mathbb{C}$  such that \begin{eqnarray}(Id \otimes \Delta )\circ \Delta= (\Delta \otimes Id)\circ \Delta \nonumber \\   (Id\otimes \varepsilon)\circ \Delta =Id= (\varepsilon \otimes Id)\circ \Delta \nonumber \end{eqnarray}

Let $C=(\mathbb{C}[x], \Delta, \varepsilon)$ be the  coalgebra of complex polynomials  with the divided power  structure
$\Delta(x^{n})=\sum_{i=0}^{n} \binom{n}{i} x^{n-i}\otimes x^{i},\;\;\varepsilon(1)=1,\;\varepsilon(x^{n})=0$ for $n>0$. The formulation of this structure comes from both  the expansion  $(x+y)^{n}=\sum_{k=1}^{n} \binom{n}{k} x^{k}y^{n-k}$ and  formula $B_{n}(x+y)=\sum_{k=1}^{n}  \binom{n}{k} B_{k}(x)B_{n-k}(y) $  where $B_{k}$'s are Bell polynomials introduced by E.T. Bell in \cite{BELL}.
  Assume  that $T$ is the operator of differentiation on $\mathbb{C}[x]$. Then $T$ satisfies
\begin{equation}\label{FQ}
(T\otimes T)\circ \Delta=\Delta \circ T^{2}
\end{equation}
So we have the  following  commutative  diagram;
\[ \begin{CD}
C \otimes C @< \Delta <<C           \\
@A {T\otimes T}AA  @AA T^{2} A  \\
C \otimes C @<\Delta <<  C
\end{CD}
\]

Now assume that $A$  is  a (not necessarily unital) algebra with multiplication $m: A\otimes A \rightarrow A$ and $T$ is  a linear map on $A$. By reversing the direction of arrows in the above diagram and  replacing the coproduct  $\Delta$ of $C$  by product m of $A$, we find the  following commutative  diagram:
\[ \begin{CD}
A \otimes A @> m >>A           \\
@V{T\otimes T}VV  @VV T^{2} V  \\
A \otimes A @>m >>  A
\end{CD}
\]

This means  that $T$  is a linear map on algebra $A$  which satisfies:
\begin{equation}\label{TP}
T(x)T(y)=T^{2}(xy)
\end{equation}
In fact,  motivated by the classical operator of differentiation, we construct (\ref{FQ})   as a coalgebraic functional equation on an arbitrary coalgebra.
This  equation naturally gives us the  equation (\ref{TP}), as  an algebraic  functional equation for  linear maps on a complex algebra $A$. If we wish to consider  (\ref{TP}) on a $C^{*}$  algebra $A$, it is natural to add the  symmetric condition $T(x^{*})=(T(x))^{*}.$\\

In this  paper we are mainly interested  in the functional equation

\begin{equation}\label{EQ}
\begin{aligned}
T(x)T(y)=T^{2}(xy) \\
T(x^{*})=(T(x))^{*}
\end{aligned}
\end{equation}

where $T$ is  a (not necessarily continuous)  linear map on a $C^{*}$  algebra $A$  and $T^{2}=T\circ T$. An operator which satisfies  (\ref{EQ}) is  called  a \textit{partial multiplier}.\\
We  observe that a partial multiplier is  automatically  continuous.
 Despite of its pure algebraic nature, we will see that for certain $C^{*}$ algebras, this  functional equation has a geometric interpretation in terms of inner product preserving maps, see Proposition 2.\\

 Our reason that we choose the name "partial multiplier" for such operators is that an injective partial multiplier on an algebra $A$ can be considered as an element of multiplier algebra of $A$,  see (\ref{inj})  Proposition 1. Another reason for this name is that a partial multiplier $T$ on a unial $C^{*}$ algebra $A$ is equal to multiplication by $T(1)$, provided that we restrict $T$ to $T(A)$. \\

 Obviously for every $\lambda \in \mathbb{R}$, the trivial linear map $T(x)=\lambda x$ is  a partial multiplier. In this paper we are
interested in conditions on a $C^{*}$  algebra $A$, under which every partial multiplier  is  necessarily  a trivial map.
\begin{thm1}
Every partial  multiplier on a $C^{*}$  algebra $A$ is trivial if $A$  satisfies  any one of the following conditions:
\begin{enumerate}[(I)]
\item $A$  is a simple  $C^{*}$-algebra. \label{S}\\

 \item $A$ is unital, with trivial center, and has  a faithful trace such that  each zero-trace element lies in the  closure of the span of  commutator  elements. \label{SS}\\

\item $A=B(H)$  where H is a  separable Hilbert  space. \label{SSS}
\end{enumerate}
\end{thm1}
\
\\
Our next Theorem, is a characterization of  all formally real fields $F$, in terms of the functional  equation (\ref{SE}) on the matrix algebra $M_{n}(F)$:
\begin{thm2}
A field $F$ is  a formally real field if  and only if for each $n\in \mathbb{N}$ the equation (\ref{SE}) has only  trivial solution $T(x)=\lambda x$ for some $\lambda \in F$.
\end{thm2}

\section{\bf Preliminaries}

In this section we give some  definitions and notations.  For a  $C^{*}$  algebra $A$, a positive  linear  map on $A$ is  a linear  map $T$  such that  $T(x)\geq 0$ for $x\geq 0$. A faithful (positive) trace on  $A$ is  a  bounded linear  map $tr: A\rightarrow \mathbb{C}$  such that $tr(xy)=tr(yx)$ and $tr(x)>0$ for $x > 0$ ($tr(x)\geq 0\; \text{for}\; x\geq 0$). A zero trace element is an element $x\in A$ with $tr(x)=0$. Elements of the  form $xy-yx$ are  called  commutators. For  a $C^{*}$  algebra $A$ with a faithful trace $tr$ we define an inner product $<.>_{tr}$ on $A$ with $<a,b>_{tr}=tr(ab^{*})$.\\

 The  dual and  bidual of  $A$ is  denoted by $A^{*}$ and $A^{**}$, respectively. A bounded linear map $T$ on $A$ induces  natural linear maps $T^{*}$ and $T^{**}$ on $A^{*}$ and $A^{**}$, respectively. The  space $A^{**}$ is  a  $C^{*}$  algebra  with the Arens  product  and  a natural involution.
The  Arens  product is  defined in three stages as follows: (for  more information on Arens  product see \cite{CY})
\begin{itemize}
\item For $f\in A^{*}$ and $x\in A$ define $<f,\;x> \in A^{*}$ with $<f,\;x>(y)=f(xy)$ for  $y\in A$.
\item For $F\in A^{**}$  and $g \in A^{*}$ define $[F,\;g]\in A^{*}$ with \newline $[F,\;g](x)=F(<g,\;x>)$ for $x\in A$.
\item For $F,G \in A^{**}$ the Arens product  $F.G\in A^{**}$ is defined   with \newline $F.G(f)=F([G,\;f])$ for $f\in A^{*}$.
\end{itemize}
The involution on $A^{**}$ is defined  as  follows:  $F^{*}(\phi)=\overline{F(\phi^{*})}$,   $\phi^{*}(a)=\overline{\phi(a^{*})}$, where  $F\in A^{**},\;\phi\in A^{*},\;a\in A$. Note that the natural involution is available when $A$ is  a  $C^{*}$ algebra but not in general when $A$ is  just  a  Banach $*$ algebra. Then $A^{**}$ is  a unital  $C^{*}$ algebra which contains $A$ as a $C^{*}$ subalgebra, via the natural imbedding of $A$ into $A^{**}$. The multiplier algebra of $A$, denoted by $\mathcal{M}(A)$, is  the idealizor of $A$ in $A^{**}$, that is the  algebra $\{z\in A^{**} \mid zA\subseteq A\;\; \& \;Az\subseteq A \}$.\\

 A pair $(L,\;R)$ of linear maps on $A$ is called a double centralizer if $R(x)y=xL(y)$ for $x,y \in A$. The  space of double centralizers on $A$ is  a $C^{*}$  algebra with  natural  operations and is isomorphic to the multiplier algebra $\mathcal{M}(A)$. In fact for every double centralizer $(L,\;R)$ on $A$ there is a unique element $a\in \mathcal{M}(A)$ such that $L(x)=ax,\;R(x)=xa$ for $x\in A$. For the algebra of double centralizers of a $C^{*}$ algebras see \cite{CENTRAL}.

\section{\bf Partial multipliers}

Let  $A$ be  a  $C^{*}$  algebra. A  linear  map $T:A\rightarrow A$ is  called a partial multiplier if $T$ satisfies  (\ref{EQ}).  Some  algebraic properties of partial multipliers  are  as follows:

\begin{prop1} \label{GEO}
Let $T$ be  a partial  multiplier on a  $C^{*}$ algebra $A$. Then
\begin{enumerate}[(a)]

\item  T is a bounded operator,  and $T^{**}$ is  a partial multiplier on $A^{**}$. \label{bo}\\

\item  $\ker(T)$ is  a closed two sided ideal in $A$. \label{ide}\\

\item    $\prod_{i=1}^{n} T(x_{i})=T^{n}(\prod_{i=1}^{n}x_{i})\; \text{where}\; x_{i}\in A \;\text{for}\; i=1,2,\ldots,n.$\\ \label{uab}

\item If $T$ is an injective operator then $(T,\;T)$ is  a  double  centralizer on $A$. \label{inj}
\end{enumerate}
\end{prop1}
\begin{proof}
 To prove (\ref{bo}) assume that $T$ is  a partial multiplier on $A$. Then $T^{2}(xx^{*})=T(x)(T(x))^{*}$ so $T^{2}$ is  a positive map. Since $T^{2}$ is a positive map on a $C^{*}$  algebra, it is  a bounded operator, see '\cite[page260]{DAVID}'. This implies that $T$ is a bounded operator too because \begin{equation*}\parallel T(x) \parallel^{2}=\parallel T(x)(T(x))^{*}\parallel =\parallel T^{2}(xx^{*})\parallel \leq \parallel T^{2}\parallel \parallel xx^{*}\parallel= \parallel T^{2}\parallel \parallel x\parallel^{2}.\end{equation*}
We  omit the  proof  of  the  last  part of (\ref{bo}) since it is  a  mimic  of the proof of Theorem 6.1 in  \cite{CY}

Now we prove (\ref{ide}). Assume that $T(x)=0$. Then for each $y\in A$ \begin{center}$T(xy)(T(xy))^{*}=T(xy)T(y^{*}x^{*})=T^{2}(xyy^{*}x^{*})=T(x)T(yy^{*}x^{*})=0$\end{center}
This shows that $\ker T$ is a right ideal in $A$.
On the other hand $\ker T$ is a $*$-subspace of $A$, since $T$ is a symmetric operator.
This  shows that $\ker T$ is  a two sided ideal in $A$.

To  prove the remaining parts of the Proposition, without loss of generality, we may assume that $A$ is a unital algebra. Otherwise we consider  the extension  $T^{**}$ of $T$ on $A$ as a linear operator on  $A^{**}$. So in the remaining part of the proof, in the non unital case, $T$ is implicitly used for $T^{**}$. This implicit usage of $T$ for $T^{**}$ is legal since the restriction of $T^{**}$ to $A$ is equal to $T$.

To prove (\ref{uab})  we first note that $T(1)^{k}T(x)=T^{k+1}(x)$ for  all $k\in \mathbb{N}$. Now we prove (\ref{uab})
by induction on $n$. Assume that the statement is true for  all $k\leq n-1$. Then \begin{multline*}
\prod_{i=1}^{n} T(x_{i})=T^{n-1}(\prod_{i=1}^{n-1} x_{i})T(x_{n})=\\
T(1)^{n-2}T(\prod_{i=1}^{n-1} x_{i})T(x_{n})=T(1)^{n-2}T^{2}(\prod_{i=1}^{n} x_{i})=\\
T(1)^{n-2}T(T(\prod_{i=1}^{n} x_{i}))=T^{n-1}(T(\prod_{i=1}^{n} x_{i})=T^{n}(\prod_{i=1}^{n}x_{i}).\end{multline*}

To prove (\ref{inj}) we note  that for all $x,y\in A$ \begin{equation}\label{TTTTT} T(x)T^2(y)=T^2(x)T(y) \end{equation}
since each side of the equality is equal to $T(x)T(1)T(y)$. In the non unital case, this $T(1)$ is replaced by $T^{**}(1)$. Now (\ref{TTTTT}) implies that  \begin{equation}\label{air}T^2(xT(y))=T^2(T(x)y)\end{equation} Since $T$ is injective we conclude that $xT(y)=T(x)y$. Then $(T, T)$ is  a double centralizer on $A$.
\end{proof}

In the  next result we give a geometric interpretation for partial multipliers. First we need
 the following lemma which is  proved in \cite[Theorem 1]{CH}:
\begin{lem1}
 Assume  that $T$ is a linear map on a complex inner product space $V$ which preserves orthogonality. Then there is  a real number $k$ such that\newline $<T(x), T(y)>=k <x, y>$ for  all $x,y\in V$.
\end{lem1}

In the  following  result  $<x, y>_{tr}=tr(xy^{*})$, is the inner product  induced from a faithful trace.
\begin{pro2}
Assume that a $C^{*}$ algebra $A$  has a faithful trace such that every zero trace element lies in the closure of span of commutator elements. Let $T$ be  a partial multiplier on $A$. Then there is  a $\lambda \in \mathbb{C}$  such that  $<T(x), T(y)>_{tr}=\lambda <x,y>_{tr}$, for all $x,y\in A$
\end{pro2}
\begin{proof}
  By above Lemma, it is sufficient to prove that $T$ preserves the orthogonality with respect to the inner product $<.,.>_{tr}$.  Note that for every $a,b\in A$,  $tr(T^2(ab-ba))=tr(T(a)T(b)-T(b)T(a))=0$.  Hence the functional $tr\circ T^2$ vanishes  on the closure of the span of commutator elements. Next suppose that $tr(xy^{*})=0$ for some  $x,y \in A$. Then $tr(T^2(xy^{*}))=0$. Hence $tr(T(x){T(y)}^{*})=0$. So $T$ preserves the orthogonality with respect to $<.,.>_{tr}$. This completes the proof of the Proposition.

\end{proof}

\begin{NEW}
To prove (\ref{S}), assume that $T$ is a non zero partial multiplier on a simple $C^{*}$  algebra $A$. By proposition 1 (\ref{ide}),  the kernel of $T$ is  an ideal in $A$ so $\ker T=\{ 0 \}$, that is $T$ is injective. By (\ref{inj}), $(T, T)$ is  a double centralizer on $A$. This  means that there is  an element $z\in \mathcal{M}(A)$ such that $zx=xz$ for all $x\in A$. By the  following argument we conclude that $z$ is a multiple of the identity element of  $\mathcal{M}(A)$. (Communicated to us by   Professor J. Rosenberg)\\
 Every $C^{*}$ algebra has an irreducible representation on a Hilbert space $H$, see \cite[Corollary I.9.11]{DAVID}. Since $A$ is simple this  representation is injective. So $A$ is an irreducible subalgebra of $B(H)$. From an equivalent definition of the multiplier algebra which is  mentioned in \cite[Proposition 2.2.11]{WO}, we have that $\mathcal{M}(A)$ is  the idealizer of $A$ in $B(H)$.
Moreover irreducibility of $A$ in $B(H)$ implies that the centralizer of $A$ in $B(H)$ reduces to one dimensional scalars $\mathbb{C}.1$,\cite[Lemma I.9.1]{DAVID}. This  obviously shows that $z$ is  a multiple of the identity. Then the partial multiplier $T$
is in the form $T(x)=\lambda x$ for some $\lambda \in \mathbb{C}$. Since $T$ is  a symmetric operator, $\lambda$ is  a real number. This complete the proof of (\ref{S}).\\

Assume that $T$ is  a partial multiplier on a $C^{*}$ algebra $A$ which satisfies  the hypothesis of (\ref{SS}). Then $T$ is  injective by Proposition 2 and $\mathcal{M}(A)=A$ since $A$ is unital. Then,  similar to the above situation,  $(T, T)$ is a double  centralizer for $A$.  Since $\mathcal{M}(A)=A$ we have $T(x)=zx=xz$ for  a central element $z\in A$. Since $A$  has trivial center we conclude that the  symmetric  linear map $T$ is in the form $T(x)=\lambda x$ for some $\lambda \in \mathbb{R}$. This proves (\ref{SS}).\\

 The same argument as above  also shows that an injective partial multiplier on $B(H)$ is  a trivial map. Then to prove (\ref{SSS}) we assume that $\ker T$ is  non trivial, then  we will obtain  a contradiction.\\
Since  $K(H)$, the space of compact operators on $H$, is the unique closed two sided ideal in $B(H)$, we may assume that $\ker T= K(H)$. Then T induces the quotient operators
$\widetilde{T}: B(H)/K(H)\rightarrow B(H) \text{and}\;\; \widehat{T} : B(H)/K(H) \rightarrow B(H)/K(H)$. The  Calkin algebra $B(H)/K(H)$ is  a simple algebra and $\widehat{T}$ is  a partial multiplier on $B(H)/K(H)$. Then (\ref{S}) implies
that  $\widehat{T}(a)\equiv \lambda a$ for some $\lambda \in \mathbb{R}$. If $\lambda=0$ then $T(x)$ is a compact operator, for all $x\in B(H)$. Thus $T^{2}=0$ since $\ker T= K(H)$. So $T(x){T(x)}^{*}=T^{2}(xx^{*})=0$, then $T(x)=0$ for all $x\in B(H)$. Now assume that $\lambda \neq 0$. After a resealing $T:= T/ \lambda$ we can assume that  $\widehat{T}$ is the identity operator. This means that $T(x)-x \in K(H)= \ker T$ for all $x\in B(H)$. Then $T^{2}=T$ so $T$ is  a $C^{*}$ morphism on $B(H)$, since $T$ satisfies (\ref{EQ}).  Let $\pi :B(H)\rightarrow B(H)/K(H)$ be the canonical map. We have
$\pi \circ \widetilde{T}=Id$, the identity operator on the Calkin algebra. This is  a contradiction by each of the following arguments:
\begin{itemize}
\item It is well known that the Calkin algebra can not be embedded in $B(H)$, see \cite[page 41]{EV}.
\item The equality $\pi \circ \widetilde{T}=Id$ implies that the  following short exact sequence is splitting:
\begin{equation*}
0\rightarrow K(H)\rightarrow B(H) \rightarrow B(H)/K(H) \rightarrow 0
\end{equation*}
On the other hand a splitting short exact sequence of $C^{*}$ algebras, gives us a splitting  exact sequence of their $K$-theory,   see \cite[Corollary 8.2.2]{WO}. \\
In particular for $i=0,1$  we would obtain that
 $$K_{i}(B(H))\approx K_{i}(K(H))\bigoplus K_{i}((B(H)/K(H))$$ This is a contradiction, see the  catalogue of K-groups in \cite[pages 123]{WO}.

\end{itemize}
This  completes the proof of  theorem 1. \hspace{6cm}$\Box$

\end{NEW}

To prove Theorem 2 we need the following lemma:
\begin{lem2}
Assume that $T$ is  a linear operator on $M_{n}(F)$ which satisfies (\ref{SE}) then, for all $x\in \ker(T)$ and for all $y\in M_{n}(F)$, $xy$  and $yx$ belong to $\ker(T^{2})$. Moreover, for every $k\in \mathbb{N}$, $T^{k}$ satisfies (\ref{SE}).
\end{lem2}

\begin{proof}
The last part of the Lemma is obvious. We prove the first part. Assume that  $T(x)=0$. Then for each $y\in M_{n}(F)$ we have $T^{2}(xy)=T(x)T(y)=0$. Moreover $(T^{2}(yx))^{tr}=(T^{2}(yx)^{tr})=T(x^{tr})T(y^{tr})=0$, hence $T^{2}(yx)=0$.
\end{proof}

\begin{prooff}
Assume  that  $F$ is  a  formally real  field and  $T$  is  a  linear  map  on $M_{n}(F)$ which satisfies (\ref{SE}). We prove that if $T\neq 0$ then $T$ is  an injective operator. On the other hand  $T$  satisfies (\ref{air}) hence   injectivity of $T$ implies $xT(y)=T(x)y,\;\;\forall x,y \in M_{n}(F)$. Then $T(x)=xT(Id)=T(Id)x$. So $T(Id)$, being   a central element in $M_{n}(F)$, is an scalar element. This  would complete the proof of one side of the Theorem. Now  assume that $T\neq 0$ and is  not injective. Since $\ker(T^{i})_{i\in \mathbb{N}}$ is  an increasing sequence of sub vector spaces of $M_{n}(F)$,  there exists  $k\in \mathbb{N}$ such that $\ker(T^{k})=\ker (T^{2k})$.  So the above Lemma implies that $\ker T^{k}$ is  an ideal in $M_{n}(F)$. This  shows that $\ker (T^{k})=M_{n}(F)$, since $M_{n}(F)$ is a simple algebra.   Then $T$ is  a nilpotent operator. Note that  (\ref{uab}) of Proposition $2$ is applicable here. Then for all $x\in M_{n}(F)$, $T(x)$ is  a  nilpotent matrix. So the  image of $T:M_{n}(F)\to  M_{n}(F)$ is  a  subvector  space of  $M_{n}(F)$ which  consists of nilpotent matrices together with zero, moreover  it is  closed under transpose operator. This  implies that the image of $T$ is the  zero space, that is $T$ is identically zero, a contradiction to $T\neq 0$.   The reason is that for  a  formally real  field $F$, the  zero matrix is  the only nilpotent  symmetric or  anti symmetric  matrix. So the above contradiction  completes the proof of one  side of the  Theorem. Now  we prove the  converse:  Assume that $F$ is  not  a  formally real field. Then there are $b_{1},b_{2},\ldots,b_{n}$ in  $F$ with $\sum_{i=1}^{n} b_{i}^{2}=0$. Put $B=(b_{i}b_{j})_{n\times n}$. Then $B=B^{tr}\neq 0$, $B^{2}=0$. Now  define $T:M_{n}(F)\to  M_{n}(F)$  with $T(A)=trace(A) B$. $T$ is  a  non  scalar operator  which  satisfies (\ref{SE}).
\end{prooff}

\section{\bf Remarks}

\begin{remark1}
There is  a wide class of  $C^{*}$ algebras which satisfy part(\ref{SS}) of Theorem 1. Apart from the matrix algebra $M_{n}(\mathbb{C})$, for every non Abelian free group $F$, the reduced $C^{*}$ algebra $C^{*}_{red}(F)$, being a unital $C^{*}$ algebra with trivial center, posses a unique faithful trace with the property that each zero trace element is  a limit of commutator elements. This is communicated to us by Professor A. Valette.
\end{remark1}

\begin{remark2}
The  following example  shows that the  assumption "faithful  trace"  in  Theorem 1,part  (\ref{SS}) cannot be  weakened  to "positive  trace":\\
Let $\mathcal{K}$ be  the  algebra  of  compact  operators on  an infinite dimensional  separable  Hilbert  space. Assume  that $A=\mathcal{K}\bigoplus \mathbb{C}$ is  the unitization of $\mathcal{K}$. Obviously $A$  is  a unital  algebra  with trivial center. Then define $tr: A \rightarrow \mathbb{C}$
with $tr(x,\; \lambda)= \lambda$. This  is  a positive  but  not  faithful trace on  $A$. Every  zero trace element is  in the  form $(T,\;0)$ which is   a sum of three   commutator element. Because every  compact  operator on  an infinite dimensional separable Hilbert space is  a  sum of three commutators, see \cite{FON}. The  operator $T(x,\; \lambda)=(0,\; \lambda)$ is  a nontrivial  partial  multiplier on $A$. This  shows that we cannot replace the assumption "faithful trace" in (\ref{SS}) by  the  weaker  assumption "positive trace".
\end{remark2}

\begin{remark3}
As  a  consequence of part (\ref{SS}) of the  main Theorem we conclude that the existence of  a nontrivial idempotent $C^{*}$  morphism on  a  unital $C^{*}$  algebra $A$ with trivial center, is  an obstruction for  $A$ to  have  a  faithful trace  with the property that each zero trace element is a sum of commutator elements.
\end{remark3}

\begin{remark4}
We  observed in Proposition  1 that every partial multiplier on  a  $C^{*}$  algebra is automatically continuous.
The  following example shows that if we remove the  symmetric condition $T(x^{*})=(T(x))^{*}$ from the definition of partial multiplier, we may loose the automatic continuity:  Let  $H$  be an infinite  dimensional Hilbert space  and  $A=B(H\bigoplus H)$. Then  each element of $A$ is in the  form $\begin{pmatrix}
X& Y\\
Z& W
\end{pmatrix}$ where $X, Y,Z$  and $W$ are  elements  of  $B(H)$.  Assume  that $\phi$ is  an unbounded  functional on $B(H)$.\\
 Then $T(\begin{pmatrix}
X& Y\\
Z& W
\end{pmatrix})=\begin{pmatrix}
0& \phi(X)Id\\
0& 0
\end{pmatrix}$
is an unbounded operator on  $A$ which satisfies $T(x)T(y)=T^{2}(xy)$.

\end{remark4}

\section*{\bf Acknowledgments}

We are grateful to the referee for very valuable and interesting comments.




\begin{thebibliography}{20}

\bibitem{BELL} E.T.Bell, \emph {Partitions Polynomials}, Annals of Math. \textbf{29}, (1927), 38-46.

\bibitem{CH} J. Chmieli\'{n}ski, \emph{Linear  mappings  approximately preserving orthogonality}, J. Math. Anal. Appl. \textbf{304} (2005), 158-169.

\bibitem{CY} P.Civin, B.Yood,  \emph{The second  conjugate space of an algebra  as an algebra}, Pacific  Journal of  Math,  \textbf{11}, (1961). no.3, 847-870

\bibitem{DAVID} K.I. Davidson, \emph{$C^{*}$ Algebras by  Examples}, Fields  Institute  Monograph, 1996
\bibitem{EV} D.E.Evans, Y. Kawahigashi,    \emph{Quantum  Symmetries on Operator  Algebras},  Oxford  Mathematical Monographs. Oxford Science publications, The Clarendon Press, Oxford University Press, New York, 1998.
\bibitem{FON} P. Fan, C. K. Fong \emph{Which  operators are the  self  commutators of  compact operators}. Proc. Amer. Math. Soc. 80 (1980), 58-60.
\bibitem{CENTRAL} B.E.Johnson  \emph{An introduction to the theory of centralizers}. Proc. London. Math. Soc. 14 (1964), 299-320



\bibitem {WO} N.E. Wgge-Olsen, \emph{ K-theory  and $C^{*}$  algebras}, Oxford University Press, 1994.





\end{thebibliography}
\end{document}